\documentclass[12pt]{article}

\usepackage{amsmath,amssymb,amsthm}
\usepackage[matrix,arrow,curve]{xy}

\newtheorem{theorem}{Theorem}

\newtheorem{prop}{Proposition}
\newtheorem{cor}{Corollary}

\theoremstyle{definition} 
\newtheorem{definition}[theorem]{Definition}

\newtheorem{example}{Example}

\begin{document}

\title{Rigged Hilbert spaces and inductive limits}
\author{S.A. Pol'shin\thanks{polshin.s@gmail.com} }
\date{}

\maketitle

\begin{abstract}
We construct a nuclear space $\Phi$ as an inductive limit of finite-dimensional subspaces of a Hilbert space $H$ in such a way that $(\Phi,H,\Phi')$ becomes a rigged Hilbert space, thus simplifying the construction by Bellomonte and Trapani.


\end{abstract}

Rigged Hilbert spaces (RHS) introduced by Gelfand and Kostyuchenko in 1955 now play an important role in the theory of eigenfunction expansions~\cite{G4,Berez,Maurin} as well as in the theory of resonant states in quantum mechanics (see~\cite{BG78,Gadella-rev} for a review). Constructions of rigged Hilbert spaces used in physics are based on spaces of analytic functions~\cite{Gadella83,Bollini96}. Recently the use of inductive limits was proposed to construct RHS~\cite{Trapani}.
In the present note we simplify the method of~\cite{Trapani} by restricting ourselves to strict inductive limits of spectra of finite-dimensional Hilbert spaces. Also we fill the gap in the construction presented in~\cite{Nagel70}.

Recall~\cite{Floret80} that a \textit{direct spectrum} is a family $\{ H_i ; i\in \mathbb{N} \}$ of LCVS together with continuous maps $\varphi_{ij}: \ H_i \rightarrow H_j$  for all $i<j$ such that the equality $\varphi_{ij} \varphi_{jk}=\varphi_{ik}$ holds for all $i<j<k$.
 
Suppose $\{ H_i \}$ is an inflating sequence of LCVS (i.e.  $H_i \subset H_{i+1}$ for all $i$), then $\Phi=\bigcup_{n-1}^\infty H_i$ endowed with inductive topology w.r.t. canonical embeddings $\varphi_i: \ H_i \rightarrow \Phi$ is called \textit{inductive limit} of the direct spectrum $\{ H_i \}$ denoted by $\lim_\rightarrow H_i$.  Then due to well-known universality of inductive topology (\cite{Schaef}, II.6.1; \cite{Bogach}, 2.3.5) we have analogous universality property of direct limit~\cite{Floret80}.

\begin{prop}\label{lim-univ}
Let $\tilde{\Phi}$ be a LCVS and let ${\tilde{\varphi}}_i :\ H_i \rightarrow \tilde{\Phi}$ be continuous maps such that $\varphi_{i\, i+1} {\tilde{\varphi}}_i ={\tilde{\varphi}}_{i+1}$ for all $i$. Then the continuous map $\lambda: \Phi\rightarrow \tilde{\Phi}$ there exists such that the diagram
\begin{equation}\label{diag}
	\xymatrix{
\tilde{\Phi} & & & \Phi\ar[lll]_\lambda \\
\ldots \ar[r] & H_i\ar[r]_{\varphi_{i\, i+1}} \ar[ul]^{{\tilde{\varphi}}_i} \ar[urr]^{\varphi_i} & H_{i+1}  \ar[r] \ar[ull]_{{\tilde{\varphi}}_{i+1}} \ar[ur]_{\varphi_{i+1}}  & \ldots
}
\end{equation}
commutes, where $\Phi=\lim_\rightarrow H_i$.
\end{prop}
Taking $\tilde{\Phi}=\mathbb{C}$ in Prop.~\ref{lim-univ} we obtain the following corollary due to Dieudonn\'e and Schwartz.
\begin{cor}\label{D-Schw}
Any linear functional on $\lim_\rightarrow H_i$ is continuous iff its restrictions to each $H_i$ are.
\end{cor}

Suppose now that $H_i \not= H_{i+1}$ for all $i$ and the topology of $H_i$ is equivalent to those induced from $H_{i+1}$, then the inductive limit is called \textit{strict}. It is well known that the strict inductive limit is complete, Hausdorff and nuclear provided all the $H_i$ are (\cite{Schaef}, II.6.4, II.6.6, III.7.4; \cite{Bogach}, 2.6.2, 2.6.5).

\begin{definition}
Let $i: \ \Phi \rightarrow H$ be an injective homomorphism of complete nuclear space $\Phi$ into Hilbert space $H$ with dense range. Then $(\Phi,H,\Phi')$ is called \textit{Gelfand triple} (or \textit{RHS}).
\end{definition}

\begin{example}
	Consider the construction of~\cite{Nagel70} based on the embedding of the space $s$ of rapidly decreasing sequences into $\ell^2$. The range of this embedding is indeed dense  in $\ell^2$ (each finite-dimensional subspace of $\ell^2$ contains an element of $s$). It remains to show that this embedding is continuous, but it follows from the following definition of topology of $s$ by seminorms equivalent to ordinary one (\cite{Bogach}, 1.3.19):
$$q_k(x)^2=\sum_{n=1}^\infty n^{4k} |x_n|^2.$$
In particular, we can take functional realization $s \cong C^\infty (\mathbf{T})$ and $\ell^2 \cong L^2 (\mathbf{T})$. 
	\end{example}

Let $H\cong \ell^2$ be a separable infinite-dimensional Hilbert space with basis $\{e_i, i\in \mathbb{N}\}$, then for any $x=\sum_{i=1}^\infty x_i e_i \in H$ we can define its projections $P_n x=\sum_{i=1}^n x_i e_i \in H_n=\mathrm{span}\{ e_i, i=1,\ldots,n\}$. Let $\{ \alpha_i, i\in \mathbb{N}\}$ be an increasing sequence of natural numbers, we see that $\Phi=\lim_\rightarrow H_{\alpha_i}$ is dense in $H$ since $\| x-P_n x  \|_H \rightarrow_{n\rightarrow\infty} 0$. 
Then from Prop.~\ref{lim-univ} it follows that the continuous map $\lambda: \ \Phi\rightarrow H$ there exist and it is injective. Indeed, otherwise we have $\lambda \varphi_{\alpha_i} (x)=0$ for some $i$ and $x\not= 0$ in $H_{\alpha_i}$, but this violates the injectivity of $\tilde{\varphi}_i$ in the diagram~(\ref{diag}).

The above construction may be reversed in the following sense. Let $\Phi=\lim_\rightarrow H_i$ be a strict inductive limit of a sequence of \textit{finite-dimensional} Hilbert spaces. Consider the bilinear form on $\Phi$ defined by $\langle x,y\rangle_\Phi = \langle x,y\rangle_{H_i}$ for $x,y \in H_i$. This definition is correct since all the $\varphi_{ij}$ are isometries and $\langle -,- \rangle_\Phi$ is continuous due to Corollary~\ref{D-Schw}. Then we can take Hilbert space $H$ to be (isomorphic to) the completion of $\Phi$ w.r.t. this norm in the spirit of original construction presented in~\cite{G4}; see also~\cite{Floret84} for another definition of continuous norm on inductive limit spaces. So we have proved the following proposition.
\begin{prop}
Consider the following data:
\begin{itemize}
	\item Separable infinite-dimensional Hilbert space $H$.
	\item Strict inductive limit of a sequence of \textit{finite-dimensional} Hilbert spaces.
\end{itemize}
Then starting from one of these data we can construct another one in such a way that $(\Phi,H,\Phi')$ becomes a RHS.
\end{prop}

\begin{example}
Let  $H_n$ be the space of polynomials on $\xi\in\mathbf{R}$ on degree $\leq n$ and let $H=L^2 (\mathbf{R})$. Consider continuous embedding $H_n \rightarrow H: \ P(\xi) \mapsto e^{-x^2 /4} P(\xi)$, then due to the density theorem (\cite{Antosik}, Theorem 4.5.1) we see that $\Phi=\lim_\rightarrow H_n$ is dense in $H$, so $(\Phi,H,\Phi')$ is a RHS.

\end{example}

		\vskip1cm
		
	\begin{flushright} 
		Institute for Theoretical Physics, \\ NSC Kharkov Institute of Physics and Technology, \\  Akademicheskaia St. 1, 61108 Kharkov, Ukraine
 \end{flushright}

\end{document}